\documentclass[11pt]{article}
\usepackage[utf8]{inputenc}
\usepackage[english]{babel}

\usepackage[margin=2.5cm]{geometry}
\usepackage{authblk}
\usepackage{graphicx}
\usepackage{tikz}
\usepackage{wrapfig}
\usepackage{amsmath}
\usepackage{amssymb}
\usepackage{amsthm}
\usepackage{thmtools}
\usepackage{hyperref}
\usepackage{cleveref}
\hypersetup{colorlinks = true, linkbordercolor = {white}}

\usepackage{amsfonts}
\usepackage{rotating}
\usepackage[backend=bibtex,style=numeric]{biblatex}
\usepackage{setspace}
\usepackage{enumerate}
\usepackage [autostyle, english = american]{csquotes}
\MakeOuterQuote{"}
\newtheorem{theorem}{Theorem}[section]
\newtheorem{prop}{Proposition}[section]
\newtheorem{lemma}[prop]{Lemma}
\theoremstyle{remark}
\newtheorem{remark}{Remark}[section]
\usepackage{pifont}

\newcommand{\HD}{\dim_{\text{H}}}
\newcommand{\PD}{\dim_{\text{P}}}
\newcommand{\AD}{\dim_{\text{A}}}
\newcommand{\LBD}{\underline{\dim}_{\text{B}}}
\newcommand{\UBD}{\overline{\dim}_{\text{B}}}
\newcommand{\LMBD}{\underline{\dim}_{\text{MB}}}
\newcommand{\HDthm}{\dim_{\emph H}}
\newcommand{\PDthm}{\dim_{\emph P}}
\newcommand{\ADthm}{\dim_{\emph A}}
\newcommand{\LBDthm}{\underline{\dim}_{\emph B}}
\newcommand{\UBDthm}{\overline{\dim}_{\emph B}}
\newcommand{\LMBDthm}{\underline{\dim}_{\emph {MB}}}

\begin{document}

\title{COINCIDENCE AND NONCOINCIDENCE OF DIMENSIONS IN COMPACT SUBSETS OF $[0,1]$}
\author{Andrew Mitchell and Lars Olsen}
\affil{School of Mathematics and Statistics, University of St Andrews, UK}

\maketitle

\begin{spacing}{1.2}
\begin{abstract}

We show that given any six numbers $r,s,t,u,v,w \in (0,1]$ satisfying $r \leq s \leq \min(t,u) \leq \max(t,u) \leq v \leq w$, it is possible to construct a compact subset of $[0,1]$ with Hausdorff dimension equal to $r$, lower modified box dimension equal to $s$, packing dimension equal to $t$, lower box dimension equal to $u$, upper box dimension equal to $v$ and Assouad dimension equal to $w$. Moreover, the set constructed is an $r$-Hausdorff set and a $t$-packing set.

\end{abstract}

\section{Introduction}

In this paper we consider the relations between six fundamental notions of fractal dimension, and provide a general construction of a class of compact subsets of $[0,1]$ in which each of these dimensions can take any permissible value. The dimensions we consider are: the Hausdorff dimension, denoted by $\HD$; the packing dimension, denoted by $\PD$; the lower and upper box dimensions, denoted by $\LBD$ and $\UBD$, respectively; the lower modified box dimension, denoted by $\LMBD$ (we note that there is also an upper modified box dimension, but this coincides with the packing dimension); and the Assouad dimension, denoted by $\AD$. We give the exact definitions in Section \ref{dims}. For a bounded subset $X$ of $\mathbb{R}^d$, these dimensions satisfy the following chain of inequalities.

\vspace{5mm}

\begin{center}
{\begin{tikzpicture}
\node[left] at (0,0) {$\HD X \leq \LMBD X$};
\node at (0.1,0.42) [rotate = 40]{$\leq$};
\node at (0.1,-0.42)[rotate = -40]{$\leq$};
\node at (0.87,0.86) {$\LBD X$};
\node at (0.87,-0.86){$\PD X$};
\node at (1.6,0.42) [rotate = -40]{$\leq$};
\node at (1.6,-0.42) [rotate = 40]{$ \leq$};
\node[right] at (1.7,0) {$\UBD X \leq \AD X$.};
\end{tikzpicture}}
\end{center}

\noindent There are many examples in the literature of sets in which one or more of the above inequalities are strict, see, for example, \cite{falconer}, \cite{fraser}, \cite{nilsson-wingren}, \cite{pesin}, \cite{schmeling-shmerkin}, \cite{Spear}, \cite{tricot}. We particularly draw the reader's attention to the result of Spear \cite{Spear}. In \cite{Spear}, it is shown that for any four numbers $s,t,u,v \in (0,1)$ with $s < \min(t,u) \leq \max(t,u) < v$ it is possible to construct a Cantor set $Y$ with $\HD Y = s$, $\PD Y = t$, $\LBD Y = u$ and $\UBD Y = v$. In Spear's construction, countably many disjoint subintervals of $[0,1]$ are chosen from left to right, each of which contains a subset $Y_n$, for $n \in \mathbb{Z}^{+}$, with $\HD Y_n = s$ and $\PD Y_n = t$. The countable stability of the Hausdorff and packing dimensions gives that the union $Y = \cup_{n=0}^{\infty} Y_n$ has $\HD Y = s$ and $\PD Y = t$. However, the lack of countable stability of the box dimensions allows the $Y_n$ sets to be arranged in such a way that the box dimensions do not coincide with the Hausdorff and packing dimensions. 

Our result could, in some ways, be thought of as a generalisation of Spear's; however, there are several differences between Spear's result and ours. Firstly, the set we construct is not necessarily a Cantor set. Furthermore, in contrast to the result in \cite{Spear}, we also consider the lower modified box dimension and the Assouad dimension. In addition, we show that the set we construct has positive and finite Hausdorff and packing measures at its Hausdorff and packing dimensions, respectively. We now give the statement of our main result.

\begin{theorem} \label{main-thm}
Let $r,s,t,u,v,w \in (0,1]$ with

\begin{center} 
{\begin{tikzpicture}
\node[left] at (0,0) {$r \leq s$};
\node at (0.1,0.4) [rotate = 40]{$\leq$};
\node at (0.1,-0.4)[rotate = -40]{$\leq$};
\node at (0.5,0.75) {$u$};
\node at (0.5,-0.75){$t$};
\node at (0.9,0.4) [rotate = -40]{$\leq$};
\node at (0.9,-0.4) [rotate = 40]{$ \leq$};
\node[right] at (1,0) {$v \leq w$.};
\end{tikzpicture}}
\end{center}
\noindent Then there exists a compact set $X \subseteq [0,1]$ with
\begin{equation*}
    \HDthm X = r \text{,}
\end{equation*}
\begin{equation*}
    \LMBDthm X = s \text{,}
\end{equation*}
\begin{equation*}
    \PDthm X = t \text{,}
\end{equation*}
\begin{equation*}
    \LBDthm X = u \text{,}
\end{equation*}
\begin{equation*}
    \UBDthm X = v \text{,}
\end{equation*}
\begin{equation*}
    \ADthm X = w \text{.}
\end{equation*}
Moreover, $X$ is an $r$-Hausdorff set and a $t$-packing set.
\end{theorem}

\begin{remark}
An analogous result holds in higher dimensions too. This is the content of Theorem \ref{HD-thm}. We detail how the proof of Theorem \ref{main-thm} can be modified to give this result in Section \ref{HD analogues}.
\end{remark}

\begin{remark}
It is possible to adapt the construction to include the case $r=0$; however, in this case $\mathcal{H}^r (X) = \infty$. Moreover, if $t=0$ then $\mathcal{P}^t (X) = \infty$ also.
\end{remark}

\begin{remark}
In our construction, the lower dimension (the natural dual to the Assouad dimension) will always take value zero since the set we construct contains isolated points.
\end{remark}

The set $X = C \cup D \cup E \cup F$ that is constructed in the proof of Theorem \ref{main-thm} is the union of four sets, $C$, $D$, $E$ and $F$, which have dimensions as follows.
\begin{center}
\begin{tabular}{ c|c|c|c|c}
  & $C$ & $D$ & $E$ & $F$ \\ 
  \hline
 $\HD$ & $r$ & $0$ & $0$ & $0$ \\  
 $\LMBD$ & $r$ & $0$ & $0$ & $s$ \\
 $\PD$ & $t$ & $0$ & $0$ & $s$\\
 $\LBD$ & $r$ & $u$ & $0$ & $s$\\
 $\UBD$ & $t$ & $v$ & $0$ & $s$\\
 $\AD$ & $t$ & $v$ & $w$ & $s$
\end{tabular}
\end{center}
The finite stability of the Hausdorff, lower modified box, packing, upper box and Assouad dimensions gives that the union of these four sets has the desired dimension. The lower box dimension is not finitely stable in general, but the sets $C$, $D$, $E$ and $F$ are constructed in such a way that the lower box dimension is stable under their union. The proof that $X$ is an $r$-Hausdorff set and a $t$-packing set follows easily from the Hausdorff and packing measures of $C$, $D$, $E$ and $F$.

The paper is structured as follows. In Section \ref{dims} we give the exact definitions of each of the dimensions we study. Then, in Section \ref{construction} we present the construction of several classes of sets, to which the sets $C$, $D$, $E$ and $F$ constructed in the proof of Theorem \ref{main-thm} belong, and state their various dimension and measure properties. The proof of Theorem \ref{main-thm} is then presented in Section \ref{main-proof}. In Section \ref{Lemma proofs} we give the proofs of the dimension and measure properties of the classes of sets constructed in Section \ref{construction}. Finally, in Section \ref{discussion} we detail some further properties of the set $X$, and explain how the proof of Theorem \ref{main-thm} can be adapted to give an analogous result in higher dimensions.

\section{Fractal dimensions and measures}\label{dims}

In this section we give the definitions of each of the dimensions that we study. While the definitions of the Hausdorff, packing and box dimensions are well known, we have nonetheless decided to include them. There are several reasons for this: firstly, several different definitions and notations exist in the literature, so stating the definitions precisely leaves no ambiguity for the reader; secondly, the definitions of the Assouad and lower modified box dimensions are perhaps less well known, and providing the definitions of all of the dimensions we study allows easier comparison of their properties; finally, the box dimensions provide a motivation for the study of the lower modified box dimension. We first give the definitions of the Hausdorff and packing dimensions.

\subsection{Hausdorff and packing dimensions}

The Hausdorff and packing dimensions are two of the most widely used notions of dimension and are defined in terms of measures: the Hausdorff and packing measures. This measure theoretic structure yields a number of natural properties one might expect a dimension to satisfy. 

Let $X \subseteq \mathbb{R}^d$. For $\alpha \geq 0$ and $\delta > 0$ write
\begin{equation*}
\mathcal{H}^{\alpha}_\delta (X) = \inf \left\{ \sum\limits_i (\text{diam}(E_i))^{\alpha} \ \Big| \ X \subseteq \bigcup_i E_i , \ \text{diam} (E_i) \leq \delta \text{ for all $i$} \right\} \text{.}
\end{equation*}
The $\alpha$-dimensional Hausdorff measure of $X$ is defined by
\begin{equation*}
\mathcal{H}^{\alpha} (X) = \sup\limits_{\delta > 0} \mathcal{H}^{\alpha}_\delta (X) \text{.}
\end{equation*}
The Hausdorff dimension of $X$ is then defined by
\begin{equation*}
    \HD X = \inf \{ \alpha \geq 0  \colon \mathcal{H}^{\alpha} (X) = 0 \} = \sup \{ \alpha \geq 0 \colon \mathcal{H}^{\alpha} (X) = \infty \} \text{.}
\end{equation*}
If $\alpha = \HD X$, then $\mathcal{H}^{\alpha} (X)$ could be zero, infinity, or positive and finite. A set for which $0 < \mathcal{H}^{\alpha} (X) < \infty$ is called an $\alpha$-Hausdorff set.

The packing dimension is the natural dual to the Hausdorff dimension. While Hausdorff measure is defined in terms of \textit{coverings} of sets less than a given diameter, packing measure considers \textit{packings} of disjoint balls. For $\alpha \geq 0$ and $\delta > 0$, write
\begin{equation*}
\overline{\mathcal{P}}^{\alpha}_\delta(X) = \sup \left\{ \sum_i (2r_i)^{\alpha} \ \bigg| \begin{split} \ & \{ B(x_i,r_i) \}_i \text{ is a family of disjoint closed balls} \\ & \text{ with $x_i \in X$ and $r_i \leq \delta$ for all $i$ } \end{split} \right\} \text{.}
\end{equation*}
The $\alpha$-dimensional packing premeasure is defined by
\begin{equation*}
\overline{\mathcal{P}}^{\alpha} (X) = \inf\limits_{\delta > 0} \overline{\mathcal{P}}^{\alpha}_\delta (X) \text{.}
\end{equation*}
At this stage we notice a difference to the Hausdorff measure: $\overline{\mathcal{P}}^{\alpha} (X)$ is not $\sigma$-additive, and therefore not a measure. The packing measure is constructed from the packing premeasure by Munroe's Method I (see \cite{munroe} for more details). Namely, the $\alpha$-dimensional packing measure of $X$ is defined by
\begin{equation*}
\mathcal{P}^{\alpha} (X) = \inf \left \{ \sum\limits_i \overline{\mathcal{P}}^{\alpha} (E_i) \ \Big| \ X \subseteq \bigcup\limits_i E_i \right \} \text{.}
\end{equation*}
The packing dimension of $X$ is then defined by
\begin{equation*}
    \PD X =  \inf \left \{ \alpha \colon \mathcal{P}^{\alpha} (X) = 0 \right \} = \sup \left \{ \alpha \colon \mathcal{P}^{\alpha} (X) = \infty \right \} \text{.}
\end{equation*}
In an analogous manner to the Hausdorff dimension, if $\alpha = \PD X$ then $\mathcal{P}^{\alpha} (X)$ could be zero, infinity, or positive and finite, and a set for which the latter is true is called an $\alpha$-packing set.

The reader is referred to Falconer's textbook \cite{falconer} for a thorough discussion of the Hausdorff and packing measures and dimensions.

\subsection{Box dimensions and modified box dimensions}

Another two widely used notions of dimension are those of the lower and upper box dimensions. Contrary to the definitions of the Hausdorff and packing dimensions, they are defined by considering coverings by sets of \textit{equal} diameter. For a bounded subset $X$ of $\mathbb{R}^d$, we denote by $N_{\delta} (X)$ the least number of sets of diameter $\delta$ required to cover $X$. The lower and upper box dimensions of $X$ are then defined by
\begin{equation*}
\LBD X = \liminf\limits_{\delta \rightarrow 0} \frac{\log N_{\delta}(X)}{-\log\delta} \, \text{,}
\end{equation*}
\begin{equation*}
\UBD X = \limsup\limits_{\delta \rightarrow 0} \frac{\log N_{\delta}(X)}{-\log\delta} \, \text{.}
\end{equation*}

Box dimensions are not countably stable in general. In fact, the lower box dimension is not even finitely stable in general. However, it is possible to obtain countably stable notions of dimension through a slight modification of their definitions. This gives rise to the modified box dimensions.\\
The lower and upper modified box dimensions of a subset $X$ of $\mathbb{R}^d$ are defined by
\begin{equation*}
\LMBD X = \inf\left\{ \sup\limits_i  \LBD E_i \ \Big| \ X \subseteq \bigcup\limits_i E_i  \right\} \text{,}
\end{equation*}
\begin{equation*}
\overline{\dim}_{\text{MB}} X = \inf\left\{ \sup\limits_i  \UBD E_i \ \Big| \ X \subseteq \bigcup\limits_i E_i  \right\} \text{.}
\end{equation*}
It is well known that the upper modified box dimension coincides with the packing dimension (a proof can be found in \cite{falconer}).

\subsection{Assouad dimension}

The final notion of dimension that we study is the Assouad dimension. It was introduced by Assouad in the 1970s \cite{assouad} and has received an increasing amount of attention in the literature in recent years, for example in \cite{fraser}, \cite{fraser-troscheit}.

The Assouad dimension of a subset $X$ of $\mathbb{R}^d$ is defined by
\begin{equation*}
\begin{split}
    \AD X = \inf \Bigg\{ \alpha \colon & \text{there exist constants } c,\rho > 0 \text{ such that,}\\ & \text{for all } 0 < r < R \leq \rho \text{, we have } \sup_{x \in X} N_r \big(B (x,R) \cap X \big) \leq c \left( \frac{R}{r} \right)^{\alpha} \Bigg\} \, \text{;}
\end{split}
\end{equation*}
recall that if $A \subseteq \mathbb{R}^d$ then $N_r (A)$ denotes the least number of sets of diameter $r$ required to cover the set $A$.

The Assouad dimension is not countably stable in general. In our construction we exploit this lack of countable stability to control the Assouad and box dimensions independently of the Hausdorff, packing and lower modified box dimensions.

\section{The construction of the sets $C$, $D$, $E$ and $F$}\label{construction}

The set $X = C \cup D \cup E \cup F$ that is constructed in the proof of Theorem \ref{main-thm} is the union of four sets, $C$, $D$, $E$ and $F$. In this section we present the constructions of several general classes of sets, to which the sets $C$, $D$, $E$ and $F$ belong, and state their various dimension and measure properties; this is the contents of Lemmas \ref{C prop}-\ref{LBD lemma}. In Section \ref{main-proof} we apply Lemmas \ref{C prop}-\ref{LBD lemma} to prove Theorem \ref{main-thm}, and then in Section \ref{Lemma proofs} we give the proofs of Lemmas \ref{C prop}-\ref{LBD lemma}.

\subsection{The set $C = C (\beta,\gamma,\mathbf{n})$}

The first of our constructions is a central Cantor set. While the dimension and measure properties of central Cantor sets are well known (see, for example, \cite{cabrelli}, \cite{garcia}), we have decided to include the full construction and proof since in this particular construction we introduce notation that we regularly refer back to in the other constructions. Furthermore, we utilise properties of this construction to simplify later proofs.

Let $0 < \beta \leq \gamma \leq 1$ and let $\mathbf{n} = (n_k)_k$ be a strictly increasing sequence of positive integers. We construct the set $C (\beta, \gamma, \mathbf{n})$ as follows. Let $(c_n (\beta, \gamma, \mathbf{n}) )_n$ be the sequence defined by
\begin{equation} \label{c_n}
    c_n (\beta,\gamma,\mathbf{n}) = 
    \begin{cases}
        2^{-1/\beta} \text{ if $n_{2k} \leq n < n_{2k+1}$ for some $k \in \mathbb{Z}^{+}$}\\
        2^{-1/\gamma} \text{ if $n_{2k+1} \leq n < n_{2(k+1)}$ for some $k \in \mathbb{Z}^{+}$}
    \end{cases} \text{.}
\end{equation}
Now, let $\Sigma_2 = \{ 0 , 1 \}$, and for each $v \in \Sigma_2$ and $n \in \mathbb{Z}^{+}$ define $S_{n,v}^{\beta,\gamma, \mathbf{n}} \colon [0,1] \rightarrow [0,1]$ by
\begin{equation} \label{IFS}
    S_{n,v}^{\beta,\gamma,\mathbf{n}} (x) = c_n (\beta,\gamma,\mathbf{n}) \, x + v \left(1 - c_n (\beta,\gamma,\mathbf{n}) \right) \text{.}
\end{equation}
Next, for each finite binary word $w = w_1 \dots w_n$ write
\begin{equation} \label{I_w}
   I_w (\beta,\gamma,\mathbf{n}) = S_{1,w_1}^{\beta,\gamma,\mathbf{n}} \dots \, S_{n,w_n}^{\beta,\gamma,\mathbf{n}} ([0,1]) 
\end{equation}
and for each $n \in \mathbb{Z}^{+}$ write $\Sigma^n_2 = \left\{ w = w_1 \dots w_n \colon w_i \in \Sigma_2 \right\}$. Let
\begin{equation} \label{C_n}
    C_n (\beta,\gamma,\mathbf{n}) = \bigcup_{w \in \Sigma^n_2} I_{w} (\beta,\gamma,\mathbf{n})
\end{equation}
and finally, define
\begin{equation} \label{C}
    C(\beta,\gamma,\mathbf{n}) = \bigcap_{n=0}^{\infty} C_n (\beta,\gamma,\mathbf{n}) \text{.}
\end{equation}
It is clear that the set $C (\beta,\gamma,\mathbf{n})$ is compact. The following lemma states the key dimension and measure properties of the set $C = C (\beta,\gamma,\mathbf{n})$.

\begin{lemma} \label{C prop}
Let $0 < \beta \leq \gamma \leq 1$ and let $\mathbf{n} = (n_k)_k$ be a strictly increasing sequence of positive integers with $\frac{n_k}{n_{k+1}} \rightarrow 0$ as $k \rightarrow \infty$. Then the set $C (\beta,\gamma,\mathbf{n})$ constructed in (\ref{C}) satisfies
\begin{equation*}
    \HDthm C (\beta,\gamma,\mathbf{n}) = \LMBDthm C (\beta,\gamma,\mathbf{n}) = \LBDthm C (\beta,\gamma,\mathbf{n}) = \beta \text{,}
\end{equation*}
\begin{equation*}
    \PDthm C (\beta,\gamma,\mathbf{n}) = \UBDthm C (\beta,\gamma,\mathbf{n}) = \ADthm C (\beta,\gamma,\mathbf{n}) = \gamma \text{,}
\end{equation*}
and
\begin{equation*}
    0 < \mathcal{H}^{\beta} (C (\beta,\gamma,\mathbf{n})) < \infty \text{,}
\end{equation*}
\begin{equation*}
    0 < \mathcal{P}^{\gamma} (C (\beta,\gamma,\mathbf{n})) < \infty \text{.}
\end{equation*}
\end{lemma}

\begin{remark}
If the sequence $\mathbf{n} = (n_k)_k$ does not satisfy $\frac{n_k}{n_{k+1}} \rightarrow 0$, then the set constructed could have different dimensions. For example, if $\mathbf{n} = (n_k)_k$ is the sequence defined by $n_k = 2^k$, then one can show that the set $C (\beta,\gamma,\mathbf{n})$ has $\HD C (\beta,\gamma,\mathbf{n}) = \LMBD C (\beta,\gamma,\mathbf{n}) = \LBD C (\beta,\gamma,\mathbf{n}) = \frac{3 \log 2}{2 \log(1/\beta) + \log(1/\gamma)}$ and $\PD C (\beta,\gamma,\mathbf{n}) = \UBD C (\beta,\gamma,\mathbf{n}) = \AD C (\beta,\gamma,\mathbf{n}) = \frac{3 \log 2}{\log(1/\beta) + 2 \log(1/\gamma)}$.
\end{remark}

\subsection{The set $D = D (\beta,\gamma,\mathbf{n},\mathbf{k})$}

Our next construction could be thought of as a countable analogue of the set $C (\beta,\gamma,\mathbf{n})$ constructed in (\ref{C}). For each $n \in \mathbb{Z}^{+}$, we construct a set $D_n (\beta,\gamma,\mathbf{n},\mathbf{k})$ by taking the union of a subset of the endpoints of the component intervals of $C(\beta,\gamma,\mathbf{n})$, and then take $D (\beta,\gamma,\mathbf{n},\mathbf{k})$ to be the union of all the $D_n$ sets. The points that are included in each $D_n (\beta,\gamma,\mathbf{n},\mathbf{k})$ are determined by a sequence $\mathbf{k} = (k_n)_n$ of positive integers.

Let $\Sigma_3 = \{ 0, 1, * \}$. For each $n \in \mathbb{Z}^{+}$, let $S_{n,v}^{\beta,\gamma,\mathbf{n}}$ be as defined in (\ref{IFS}) if $v \in \{ 0 , 1 \}$, and let $S_{n,*}^{\beta,\gamma,\mathbf{n}} \colon [0,1] \rightarrow [0,1]$ denote the identity map. For each finite word $w = w_1 \dots w_n$ with letters in $\Sigma_3$, write
\begin{equation} \label{P_w}
P_w (\beta,\gamma,\mathbf{n}) = S_{1,w_1}^{\beta,\gamma,\mathbf{n}} \dots \, S_{n,w_n}^{\beta,\gamma,\mathbf{n}} (\{ 0,1 \}) \text{.}
\end{equation}
Next, for each $n \in \mathbb{Z}^{+}$, define
\begin{equation*}
    \Omega_n (\mathbf{k})= \left\{ w = w_1 \dots w_n \ \bigg| \ 
    \begin{split}
        & w_i \in \Sigma_3 \text{ for all $1 \leq i \leq n$}\\
        & \text{if $w_i = 1$ and $k_i < n$ then $w_j = *$ for all $j > k_i$}
    \end{split}
    \right\} \text{.}
\end{equation*}
Set
\begin{equation} \label{D_n}
    D_n (\beta,\gamma,\mathbf{n},\mathbf{k}) = \bigcup_{w \in \Omega_n (\mathbf{k})} P_w (\beta,\gamma,\mathbf{n})
\end{equation}
and finally, define
\begin{equation} \label{D}
    D (\beta,\gamma,\mathbf{n},\mathbf{k}) = \bigcup_{n=0}^{\infty} D_n (\beta,\gamma,\mathbf{n},\mathbf{k}) \text{.}
\end{equation}
There are only finitely many points in $D (\beta,\gamma,\mathbf{n},\mathbf{k})$ outwith any ball centred at the origin, so it follows that the set $D (\beta,\gamma,\mathbf{n},\mathbf{k})$ is compact.

\begin{lemma} \label{D prop}
Let $0 < \beta \leq \gamma \leq 1$, let $\mathbf{n} = (n_k)_k$ be a strictly increasing sequence of positive integers with $\frac{n_k}{n_{k+1}} \rightarrow 0$ as $k \rightarrow \infty$ and let $\mathbf{k} = (k_n)_n$ be a strictly increasing sequence of positive integers with $\frac{n}{k_n} \rightarrow 0$ as $n \rightarrow \infty$. Then the set $D (\beta,\gamma,\mathbf{n},\mathbf{k})$ constructed in (\ref{D}) satisfies
\begin{equation*}
    \HDthm D (\beta,\gamma,\mathbf{n},\mathbf{k}) = \LMBDthm D (\beta,\gamma,\mathbf{n},\mathbf{k}) = \PDthm D (\beta,\gamma,\mathbf{n},\mathbf{k}) = 0 \text{,}
\end{equation*}
\begin{equation*}
    \LBDthm D (\beta,\gamma,\mathbf{n},\mathbf{k}) = \beta \text{,}
\end{equation*}
\begin{equation*}
    \UBDthm D (\beta,\gamma,\mathbf{n},\mathbf{k}) = \ADthm D (\beta,\gamma,\mathbf{n},\mathbf{k}) = \gamma \text{.}
\end{equation*}
\end{lemma}
\begin{remark}
If we removed the dependence on the sequence $\mathbf{k} = (k_n)_n$ and simply took the union of all the endpoints of intervals in $C_n (\beta,\gamma,\mathbf{n})$, then the constructed set would have the same dimensions, but would not be compact.
\end{remark}
\noindent Lemma \ref{D prop} says that the construction of $D (\beta,\gamma,\mathbf{n},\mathbf{k})$ can be used to control the lower and upper box dimensions independently of the countably stable notions of dimension. In our next construction we show that by changing the sequence $\mathbf{k}$ it is possible to construct sets with different dimensions. 

\subsection{The set $E = E(\gamma,\mathbf{j})$}
By changing the conditions on the sequence $\mathbf{k}$ in $D (\beta,\gamma,\mathbf{n},\mathbf{k})$ it is possible to construct a set for which the Assouad dimension takes value $\gamma$, but all other dimensions are zero. Note that if $\beta = \gamma$, then the set $D_n (\gamma,\gamma,\mathbf{n},\mathbf{k})$ does not depend on the sequence $\mathbf{n}$, therefore the set $E = E(\gamma,\mathbf{k})$ that we construct does not have a dependence on $\mathbf{n}$. For this reason we can suppress the dependence on $\mathbf{n}$ in the definition below. For each $n \in \mathbb{Z}^{+}$, write
\begin{equation} \label{E_n}
    E_n (\gamma,\mathbf{k}) = D_n (\gamma,\gamma,\mathbf{n},\mathbf{k})
\end{equation}
and set
\begin{equation} \label{E}
    E (\gamma,\mathbf{k}) = \bigcup_{n=0}^{\infty} E_n (\gamma,\mathbf{k}) \text{.}
\end{equation}

\begin{lemma} \label{E prop}
Let $0 < \gamma \leq 1$ and let $\mathbf{j} = (j_n)_n$ be a nondecreasing sequence of positive integers with
\begin{enumerate}[(i)]
    \item $\frac{j_{n}}{n} \rightarrow 1$ as $n \rightarrow \infty$,
    \item $j_n - n \rightarrow \infty$ as $n \rightarrow \infty$.
\end{enumerate}
(For example, $\mathbf{j} = (j_n)_n$ could be the sequence defined by $j_n = \big\lfloor n^{1 + \frac{1}{n}} \big\rfloor$ for all $n$.)

\vspace{0.5em}

\noindent Then the set $E(\gamma,\mathbf{j})$ defined in (\ref{E}) satisfies
\begin{equation*}
    \HDthm E(\gamma,\mathbf{j}) = \LMBDthm E(\gamma,\mathbf{j}) = \PDthm E(\gamma,\mathbf{j}) = \LBDthm E(\gamma,\mathbf{j}) = \UBDthm E(\gamma,\mathbf{j}) = 0 \text{,}
\end{equation*}
\begin{equation*}
    \ADthm E(\gamma,\mathbf{j}) = \gamma \text{.}
\end{equation*}
\end{lemma}

\subsection{The set $F = F (\gamma,\mathbf{a},\mathbf{b})$}
Finally, we present the construction of the set from which we control the lower modified box dimension. One could think of this as a modification of a construction of Tricot \cite{tricot}. In Tricot's construction a subset of real numbers is constructed with Hausdorff dimension equal to zero and lower modified box dimension equal to one. We adapt this construction so that the lower modified box dimension takes the value $\gamma$.

We define an enumeration of closed intervals $\{ V_w (\gamma,\mathbf{a},\mathbf{b}) \}_{w \in \Sigma_2^{*}}$ indexed by finite binary words, which depend on two sequences $\mathbf{a} = (a_k)_k$ and $\mathbf{b} = (b_k)_k$. We give the exact definitions later, but one should note that $V_{\varepsilon} (\gamma,\mathbf{a},\mathbf{b})= [0,1]$, any interval $V_{w_1 \dots w_n} (\gamma,\mathbf{a},\mathbf{b})$ is contained in $V_{w_1 \dots w_m} (\gamma,\mathbf{a},\mathbf{b})$ for all $m \leq n$, and any two intervals indexed by words of the same length intersect on at most one point. For each $n$, let
\begin{equation} \label{F_k}
    F_n (\gamma,\mathbf{a},\mathbf{b}) = \bigcup_{w \in \Sigma_2^{n}} V_w (\gamma,\mathbf{a},\mathbf{b})
\end{equation}
and finally, define
\begin{equation} \label{F}
    F (\gamma,\mathbf{a},\mathbf{b}) = \bigcap_{n = 0}^{\infty} F_n (\gamma,\mathbf{a},\mathbf{b}) \text{.}
\end{equation}

We now turn towards the exact definitions of the component intervals $ V_w (\gamma,\mathbf{a},\mathbf{b})$. We require $\mathbf{a} = (a_k)_k$ and $\mathbf{b} = (b_k)_k$ to be strictly increasing sequences of positive integers that satisfy the following conditions.
\begin{enumerate}[(i)]
    \item $a_0 = b_0 = 1$,
    \item $a_j - a_i \leq b_j - b_i$ for all $i,j \in \mathbb{Z^{+}}$ with $i \leq j$,
    \item $\frac{a_k}{b_k} \rightarrow 0$ as $k \rightarrow \infty$,
    \item $\frac{b_k}{a_{k+1}} \rightarrow 0$ as $k \rightarrow \infty$.
\end{enumerate}
(For example, $\mathbf{a} = (a_k)_k$ and $\mathbf{b} = (b_k)_k$ could be the sequences defined by $a_k = 2^{k^2}$ and $b_k = (k+1) \, 2^{k^2}$ for all $k$.)

\noindent Observe that we can define a bijection $f \colon \Sigma_2^{*} \rightarrow \mathbb{Z}^{+}$ by
\begin{equation*}
\begin{split}
    f (\varepsilon) & = 1 \\
    f (w) & = 2^{-n} + \sum_{i=1}^{n} w_i \, 2^{n-i} \, \text{ for $w = w_1 \dots w_n$,}
\end{split}
\end{equation*}
where $\Sigma_2^{*}$ denotes the set of all finite binary words. Hence, for each $w \in \Sigma_2^{*}$ there is a unique positive integer $k$ such that $f^{-1} (k) = w$. We are now in a position to give the exact definition of $\{ V_w (\gamma,\mathbf{a},\mathbf{b}) \}_w$.

Let $V_{\epsilon} (\gamma,\mathbf{a},\mathbf{b})= [0,1]$ and $\alpha (\gamma) = 2^{1/\gamma}$. For each $w = w_1 \dots w_n$ ($n \geq 1$) we define $V_w (\gamma,\mathbf{a},\mathbf{b})$ inductively as follows. If there is an $m < n$ such that $n = a_{f^{-1} (w_1 \dots w_m)}$ then $V_w (\gamma,\mathbf{a},\mathbf{b})$ has length $(\alpha (\gamma))^{-b_{f^{-1} (w_1 \dots w_m)}}$. Otherwise $V_w (\gamma,\mathbf{a},\mathbf{b})$ has length equal to $(\alpha (\gamma))^{-1}$ times the length of its parent interval, $V_{w_1 \dots w_{n-1}} (\gamma,\mathbf{a},\mathbf{b})$. In both cases, the intervals are positioned to the extreme left (if $w_n = 0$) or extreme right (if $w_n = 1$) of their parent interval.

Since every point in $F (\gamma,\mathbf{a},\mathbf{b})$ is the intersection of a sequence of nested closed intervals, it follows that $F (\gamma,\mathbf{a},\mathbf{b})$ is compact.

\begin{lemma} \label{F prop}
Let $0 < \gamma \leq 1$, and let $\mathbf{a} = (a_k)_k$ and $\mathbf{b} = (b_k)_k$ be strictly increasing sequences of positive integers satisfying conditions (i)-(iv). Then the set $F (\gamma,\mathbf{a},\mathbf{b})$ constructed in (\ref{F}) satisfies
\begin{equation*}
    \HDthm F (\gamma,\mathbf{a},\mathbf{b}) = 0 \text{,}
\end{equation*}
\begin{equation*}
    \LMBDthm F(\gamma,\mathbf{a},\mathbf{b}) = \PDthm F(\gamma,\mathbf{a},\mathbf{b}) = \LBDthm F(\gamma,\mathbf{a},\mathbf{b}) = \UBDthm F(\gamma,\mathbf{a},\mathbf{b}) = \ADthm F(\gamma,\mathbf{a},\mathbf{b}) = \gamma \text{.}
\end{equation*}
\end{lemma}

\subsection{On the lower box dimension of $C \cup D \cup E \cup F$}

Our next result states that the lower box dimension of the union $C \cup D \cup E \cup F$ is equal to the maximum of the lower box dimensions of the sets $C$, $D$, $E$ and $F$.

\begin{lemma}\label{LBD lemma}
Let $r,s,t,u,v,w \in [0,1]$ be as in the statement of Theorem \ref{main-thm}, $\mathbf{n} = (n_k)_k$ as in the statement of Lemma \ref{C prop}, $\mathbf{k} = (k_n)_n$ as in the statement of Lemma \ref{D prop}, $\mathbf{j} = (j_n)_n$ as in the statement of Lemma \ref{E prop} and $\mathbf{a} = (a_k)_k$ and $\mathbf{b} = (b_k)_k$ as in the statement of Lemma \ref{F prop}. Let $C = C(r,t,\mathbf{n})$, $D = D(u,v,\mathbf{n},\mathbf{k})$, $E = E(w,\mathbf{j})$ and $F = F(s,\mathbf{a},\mathbf{b})$, as defined in (\ref{C}), (\ref{D}), (\ref{E}) and (\ref{F}), and set $X = C \cup D \cup E \cup F$. Then $\LBDthm X = u$.
\end{lemma}

We are now in a position to give the proof of Theorem \ref{main-thm}. This is the content of the next section. We return to the proofs of Lemmas \ref{C prop}-\ref{LBD lemma} in Section \ref{Lemma proofs}.

\section{Proof of Theorem \ref{main-thm}} \label{main-proof}

In this section we give the proof of our main theorem.

\textit{Proof of Theorem \ref{main-thm}.} Let $r,s,t,u,v,w \in [0,1]$ be as in the statement of Theorem \ref{main-thm}, $\mathbf{n} = (n_k)_k$ as in the statement of Lemma \ref{C prop}, $\mathbf{k} = (k_n)_n$ as in the statement of Lemma \ref{D prop}, $\mathbf{j} = (j_n)_n$ as in the statement of Lemma \ref{E prop} and $\mathbf{a} = (a_k)_k$ and $\mathbf{b} = (b_k)_k$ as in the statement of Lemma \ref{F prop}. Let $C = C(r,t,\mathbf{n})$, $D = D(u,v,\mathbf{n},\mathbf{k})$, $E = E(w,\mathbf{j})$ and $F = F(s,\mathbf{a},\mathbf{b})$, as defined in (\ref{C}), (\ref{D}), (\ref{E}) and (\ref{F}), and set $X = C \cup D \cup E \cup F$. It follows by Lemmas \ref{C prop}, \ref{D prop}, \ref{E prop}, and \ref{F prop} that the dimensions of $C$, $D$, $E$ and $F$ are as follows.
\begin{center}
\begin{tabular}{ c|c|c|c|c}
  & $C$ & $D$ & $E$ & $F$ \\ 
  \hline
 $\HD$ & $r$ & $0$ & $0$ & $0$ \\  
 $\LMBD$ & $r$ & $0$ & $0$ & $s$ \\
 $\PD$ & $t$ & $0$ & $0$ & $s$\\
 $\LBD$ & $r$ & $u$ & $0$ & $s$\\
 $\UBD$ & $t$ & $v$ & $0$ & $s$\\
 $\AD$ & $t$ & $v$ & $w$ & $s$
\end{tabular}
\end{center}
The finite stability of the Hausdorff, lower modified box, packing, upper box and Assouad dimensions gives $\HD X = r$, $\LMBD X = s$, $\PD X = t$, $\UBD X = v$ and $\AD X = w$, and it follows by Lemma \ref{LBD lemma} that $\LBD X = u$. 

It remains to show $0 < \mathcal{H}^r (X) \leq \mathcal{P}^t (X) < \infty$. Lemma \ref{C prop} gives that $0 < \mathcal{H}^r (C) < \infty$ and $0 < \mathcal{P}^t (C) < \infty$, so since $C \subseteq X$, it follows that $\mathcal{H}^r (X) > 0$ and $\mathcal{P}^t (X) > 0$ also. On the other hand, since $\HD D = \PD D = 0$ and $\HD E = \PD E = 0$ we have $\mathcal{H}^r (D) = \mathcal{P}^t (D) = 0$ and $\mathcal{H}^r (E) = \mathcal{P}^t (E) = 0$. In addition, since $\HD F = 0$ we have $\mathcal{H}^r (F) = 0$, so $\mathcal{H}^r (X) = \mathcal{H}^r (C) < \infty$. Finally, observe that if $\mathbf{m} = (m_k)_k$ is any strictly increasing sequence of positive integers with $\frac{m_k}{m_{k+1}} \rightarrow 0$ as $k \rightarrow \infty$, then $F \subseteq C (s,s,\mathbf{m})$, and so $\mathcal{P}^t (F) \leq \mathcal{P}^t (C(s,s,\mathbf{m})) \leq \mathcal{P}^s (C(s,s,\mathbf{m})) < \infty$ by Lemma \ref{C prop}. Hence $\mathcal{P}^t (X) \leq \mathcal{P}^t (C) + \mathcal{P}^t (F) < \infty$. This completes the proof of Theorem \ref{main-thm}. \qed

\section{Proofs of Lemmas \ref{C prop}-\ref{LBD lemma}} \label{Lemma proofs}

We now give the proofs of Lemmas \ref{C prop}-\ref{LBD lemma}.

\subsection{Proof of Lemma \ref{C prop}}

Let $0 < \beta \leq \gamma \leq 1$, and let $\mathbf{n} = (n_k)_k$ be as in the statement of Lemma \ref{C prop}. To ease notation, we write $c_n = c_n (\beta,\gamma,\mathbf{n})$, $I_w = I_w (\beta,\gamma,\mathbf{n})$, $C_n = C_n (\beta,\gamma,\mathbf{n})$ and $C = C (\beta,\gamma,\mathbf{n})$ throughout.

\textit{Claim 1.} $\LBD C \leq \beta$.

\textit{Proof of Claim 1.} Observe that one interval of diameter $\delta_k = \prod_{i=n_{2k}}^{n_{2k+1}} c_i = 2^{-(n_{2k+1}-n_{2k})/\beta}$ will cover any component interval of $C_{n_{2k+1}}$. Hence $N_{\delta_k}(C) \leq 2^{n_{2k+1}}$, so
\begin{equation*}
    \frac{\log N_{\delta_k} (C)}{-\log \delta_k} \leq \frac{\beta \, n_{2k+1} \log 2}{(n_{2k+1}-n_{2k}) \log 2} \rightarrow \beta \text{;}
\end{equation*}
therefore $\LBD C \leq \beta$. \qed

\textit{Claim 2.} $\AD C \leq \gamma$.

\textit{Proof of Claim 2.} Let $x \in C$ and $0 < r \leq R < \frac{1}{2}$. Furthermore, let $m(R)$ and $n(r)$ be the unique integers such that $l_{n(r)} \leq r < l_{n(r)-1}$ and $\frac{1}{2} l_{m(R)+1} \leq R < \frac{1}{2} l_{m(R)}$, where $l_n = \prod_{i=1}^{n} c_i$. Then $B(x,R)$ will intersect at most two intervals of $C_{m(R)}$. Each contains precisely $2^{n(r)-m(R)}$ intervals of $C_{n(r)}$, so $2^{n(r)-m(R)+1}$ intervals of diameter $r$ will cover $B(x,R) \cap C$; hence
\begin{equation*}
    N_r (B(x,R) \cap C) \leq 2^{n(r)-m(R)+1} = 8 \, \left( 2^{-(n(r)-1)/\gamma} \right)^{-\gamma} \left( 2^{-(m(R)+1)/\gamma} \right)^{\gamma} \leq 8 \, \left( \frac{R}{r} \right)^{\gamma} \text{,}
\end{equation*}
so we conclude that $\AD C \leq \gamma$. \qed

\textit{Claim 3.} $0 < \mathcal{H}^{\beta} (C) < \infty$ and $0 < \mathcal{P}^{\gamma} (C) < \infty$.

\textit{Proof of Claim 3.} Define a mass distribution on $C$ as follows. For each $n \in \mathbb{Z}^{+}$ and $w \in \Sigma^n$, let $\lambda_w$ denote the Lebesgue measure restricted to the interval $I_w$, normalised such that $\lambda_w (I_w) = 1$. Next, define the probability measure $\mu_n$ by $\mu_n = \frac{1}{2^n} \sum_{w \in \Sigma^n} \lambda_w$. Then there exists a probability measure $\mu$ such that $\mu_n$ converges weakly to $\mu$.

For each $x \in C$, there is a decreasing sequence of intervals $\big( I_{w_1 \dots w_n} \big)_n$ such that $x = \cap_n I_{w_1 \dots w_n}$. By definition, each such interval carries mass $2^{-n}$ and has length $l_n = \prod_{i=1}^n c_i$. Let $r > 0$, and $n(r)$ be the integer such that $l_{n(r)} \leq r < l_{n(r)-1}$. Then the ball with centre at $x$ and radius $r$ can intersect at most three component intervals of $C_{n(r)-1}$, so $\mu (B(x,r)) \leq 3 \cdot 2^{-n(r)+1}$. Since $c_i \geq 2^{-1/\beta}$ for all $i \in \mathbb{Z}^{+}$, we have
\begin{equation*}
    \frac{\mu (B(x,r))}{r^{\beta}} \leq \frac{3 \cdot 2^{-n(r)+1}}{(l_{n(r)})^{\beta}} \leq \frac{3 \cdot 2^{-n(r)+1}}{(2^{-n(r)/\beta})^{\beta}} = 6
\end{equation*}
so it follows by the Mass Distribution Principle that $\mathcal{H}^{\beta} (C) \geq \frac{1}{6}$. 

On the other hand, if $r_k = l_{n_{2k+1}}$ then $B(x,r_k)$ will cover the interval of $C_{n_{2k+1}}$ that includes $x$, so $\mu(B(x,r_k)) \geq 2^{-n_{2k+1}}$. Since $\frac{n_k}{n_{k+1}} \rightarrow 0$ as $k \rightarrow \infty$ we have $l_{n_{2k+1}} \approx 2^{-n_{2k+1}/\beta}$ for all sufficiently large $k$, so
\begin{equation*}
    \frac{\mu (B(x,r_k))}{r_k^{\, \beta}} \geq \frac{2^{-n_{2k+1}}}{{(l_{n_{2k+1}}})^{\beta}} \approx \frac{2^{-n_{2k+1}}}{(2^{-n_{2k+1}/\beta})^{\beta}} = 1
\end{equation*}
for large enough $k$. Hence we conclude by the Mass Distribution Principle that $\mathcal{H}^{\beta} (C) \leq 2^{\beta} \leq 2$. 

One can show that $0 < \mathcal{P}^{\gamma} (C)< \infty$ by similar arguments. \qed

\textit{Proof of Lemma \ref{C prop}.} It follows immediately from Claim 3 that $\HD C = \beta$ and $\PD C = \gamma$. Combining this with Claim 1 gives that $\HD C = \LMBD C = \LBD C = \beta$, and combining with Claim 2 gives that $\HD C = \UBD C =\AD C = \gamma$. This completes the proof of Lemma \ref{C prop}. \qed

\subsection{Proof of Lemma \ref{D prop}}

Let $0 < \beta \leq \gamma \leq 1$, and let $\mathbf{n} = (n_k)_k$ and $\mathbf{k} = (k_n)_n$ be as in the statement of Lemma \ref{D prop}. To ease notation we write $c_n = c_n (\beta,\gamma,\mathbf{n})$, $D_n = D_n (\beta,\gamma,\mathbf{n},\mathbf{k})$ and $D = D (\beta,\gamma,\mathbf{n},\mathbf{k})$ throughout.

\textit{Claim 1.} $\LBD D \leq \beta$ and $\AD D \leq \gamma$.

\textit{Proof of Claim 1.} We have that $D_n \subseteq C_m (\beta,\gamma,\mathbf{n})$ for all $m \geq n$, therefore $D \subseteq C (\beta,\gamma,\mathbf{n})$ and the result follows by monotonicity and Lemma \ref{C prop}. \qed

\textit{Claim 2.} $\LBD D \geq \beta$.

\textit{Proof of Claim 2.} Let $0 < \delta < 1$ and $k (\delta)$ be the integer such that $l_{k(\delta)+1} \leq \delta < l_{k(\delta)}$, where $l_{k} = \prod_{i=1}^{k} c_i$. Furthermore, let $n (\delta)$ be the largest integer such that $k_{n(\delta)} < k (\delta)$. Any interval of diameter $\delta$ can intersect at most one point in $D_{k(\delta)}$, and there are at least $2^{k(\delta)-n(\delta)}$ points in $D_{k(\delta)}$, so $N_{\delta} (D) \geq 2^{k(\delta)-n(\delta)}$. Since $c_i \leq 2^{1/\beta}$ for all $i \in \mathbb{Z}^{+}$, it follows that
\begin{equation*}
    \frac{\log N_{\delta}(D)}{-\log \delta} \geq \frac{\beta \, (k(\delta) - n(\delta))}{k(\delta)+1} \geq \frac{\beta \, k(\delta)}{k(\delta)+1} - \frac{n(\delta)}{k_{n(\delta)}}\rightarrow \beta
\end{equation*}
since $\frac{n(\delta)}{k_{n(\delta)}} \rightarrow 0$ as $\delta \rightarrow 0$. Therefore $\LBD D \geq \beta$. \qed

\textit{Claim 3.} $\UBD D \geq \gamma$.

\textit{Proof of Claim 3.} If, for each $k \in \mathbb{Z}^{+}$, $\delta_k = l_{n_{2k}-1}$ and $m(k)$ is the largest integer such that $k_{m(k)} \leq n_{2k}$, then any interval of diameter $\delta_k$ can intersect at most one point in $D_{n_{2k}}$. There are at least $2^{n_{2k} - m(k)}$ such points, so $N_{\delta_k} (D) \geq 2^{n_{2k} - m(k)}$ and it follows that
\begin{equation*}
    \frac{\log N_{\delta_k} (D)}{-\log \delta_k} \geq \frac{n_{2k} - m(k)}{(n_{2k}-n_{2k-1})/\gamma + n_{2k-1}/\beta} \geq \frac{n_{2k}}{n_{2k}/\gamma + n_{2k-1}/\beta} - \frac{m(k)}{k_{m(k)}}\rightarrow \gamma
\end{equation*}
since $\frac{n_{2k-1}}{n_{2k}} \rightarrow 0$ and $\frac{m(k)}{k_{m(k)}} \rightarrow 0$. Hence $\UBD D \geq \gamma$. \qed

\textit{Proof of Lemma \ref{D prop}.} Combining Claims 1, 2 and 3 completes the proof of Lemma \ref{D prop}. \qed

\subsection{Proof of Lemma \ref{E prop}}

Let $0 < \gamma \leq 1$, and let $\mathbf{j} = (j_n)_n$ be as in the statement of Lemma \ref{E prop}. To ease notation we write $c_n = c_n (\beta,\gamma,\mathbf{n})$, $E_n = E_n (\gamma,\mathbf{j})$ and $E = E (\gamma,\mathbf{j})$ throughout.

\textit{Claim 1.} $\UBD E = 0$.

\textit{Proof of Claim 1.} Let $0 < \delta < 1$ and $n (\delta)$ be the integer such that $l_{n(\delta)} \leq \delta < l_{n(\delta)-1}$, where $l_{n} = \prod_{i=1}^{n} c_i$. Since $\cup_{i = j_{n(\delta)}}^{\infty} E_i \subseteq [0, l_{n(\delta)}]$, one interval of diameter $\delta$ will cover $\cup_{i=j_{n(\delta)}}^{\infty} E_i$. For each nonnegative integer $m$, $\cup_{i=j_{m}}^{j_{m+1}} E_i$ contains $2^{j_m - m}$ points, so an additional $\sum_{i=1}^{n(\delta)-1} 2^{j_i - i} \leq (n(\delta) - 1) \, 2^{j_{n(\delta)} - n(\delta)}$ intervals of diameter $\delta$ will cover the remaining points in $E$. Thus $N_{\delta} (E) \leq n(\delta) \, 2^{j_{n(\delta)}-n(\delta)}$ and it follows that
\begin{equation*}
    \frac{\log N_{\delta} (E)}{-\log \delta} \leq \frac{\log n(\delta)}{n(\delta) - 1} + \frac{j_{n(\delta)} - n(\delta)}{n (\delta) - 1} \rightarrow 0
\end{equation*}
since $\frac{j_{n(\delta)}}{n(\delta)-1} \rightarrow 1$. Hence $\UBD E = 0$. \qed

\textit{Claim 2.} $\AD E = \gamma$.

\textit{Proof of Claim 2.} For any sequence $\mathbf{m} = (m_k)_k$ with $\frac{m_k}{m_{k+1}} \rightarrow 0$ as $k \rightarrow \infty$ we have $E \subseteq C(\gamma,\gamma,\mathbf{m})$, so it follows by monotonicity and Lemma \ref{C prop} that $\AD E \leq \gamma$. It remains to show the reverse inequality. To this end, let $A>0$ and $0 < \alpha < \gamma$. For each $n \in \mathbb{Z}^{+}$, $B(0,l_{n})$ will contain $E_{j_n}$. There are $2^{j_n - n}$ points in $E_{j_n}$, and any interval of length $l_{j_n}$ can intersect at most two of them, so
\begin{equation*}
    N_{l_{j_n}} \big( B(0,l_{n} ) \cap E \big) \geq 2^{j_n-n-1} = \frac{1}{2} \left( \frac{l_{n}}{l_{j_n}} \right)^{\gamma} = \frac{1}{2} \left( \frac{l_{n}}{l_{j_n}} \right)^{\gamma - \alpha} \, \left( \frac{l_{n}}{l_{j_n}} \right)^{\alpha} \text{.}
\end{equation*}
Since $j_n - n \rightarrow \infty$ as $n \rightarrow \infty$, it follows that $\frac{l_n}{l_{j_n}} \rightarrow \infty$ also. Thus there is an $N$ such that $\left( \frac{l_{n}}{l_{j_n}} \right)^{\gamma - \alpha} > 2A$ for all $n \geq N$. Hence we conclude that $\AD E \geq \gamma$. \qed

\textit{Proof of Lemma \ref{E prop}.} Combining Claim 1 and Claim 2 completes the proof of Lemma \ref{E prop}. \qed

\subsection{Proof of Lemma \ref{F prop}}

Let $0 < \gamma \leq 1$, and let $\mathbf{a} = (a_k)_k$ and $\mathbf{b} = (b_k)_k$ be as in the statement of Lemma \ref{F prop}. To ease notation, we write $V_w = V_w (\gamma,\mathbf{a},\mathbf{b})$, $\alpha = \alpha (\gamma)$, $F_k = F_k (\gamma,\mathbf{a},\mathbf{b})$ and $F = F (\gamma,\mathbf{a},\mathbf{b})$ throughout. In addition, we denote by $p(k)$ the integer such that if $f (w_1 \dots w_n) = k$ then $f (w_1 \dots w_{n-1}) = p(k)$ (i.e. if $k$ is the integer that indexes a given interval, then $p(k)$ is the integer that indexes its parent interval).

\textit{Claim 1.} $\AD F \leq \gamma$.

\textit{Proof of Claim 1.} Observe that for any sequence $\mathbf{m} = (m_k)_k$ with $\frac{m_k}{m_{k+1}} \rightarrow 0$ as $k \rightarrow \infty$, each component interval $V_w$ is included in the interval $I_w (\gamma,\gamma,\mathbf{m})$ defined in (\ref{I_w}), so $F_n \subseteq C_n (\gamma,\gamma,\mathbf{m})$ for all $n$. Therefore $F \subseteq C(\gamma,\gamma,\mathbf{m})$, so it follows by monotonicity and Lemma \ref{C prop} that $\AD F \leq \gamma$. \qed

\textit{Claim 2.} $\LBD F \geq \gamma$.

\textit{Proof of Claim 2.} For each $k \in \mathbb{Z}^{+}$, define
\begin{equation}
    \mathcal{I}_k = \left\{ V_w \colon V_w \text{ is a component interval of } F_{a_k - 1} \text{ and } V_w \subseteq V_{f^{-1}(k)} \right\} \text{.}
\end{equation}
Let $W \in \mathcal{I}_k$. Then $W \subseteq V_{f^{-1} (k)}$, so $W \subseteq V_{f^{-1} (p(k))}$. Hence the interval of $F_{a_{p(k)}}$ that contains $W$ has length $\alpha^{b_{p(k)}}$, and it follows that $W$ has length $\alpha^{-l_k}$, where $l_k = a_k - 1 - a_{p(k)} + b_{p(k)}$.

Observe that
\begin{equation} \label{limit}
    \frac{l_{k-1}}{b_{p(k)}} = \frac{a_k}{b_{p(k)}} - \frac{1}{b_{p(k)}} - \frac{a_{p(k)}}{b_{p(k)}} + 1 \rightarrow \infty \text{,}
\end{equation}
since the first term tends to infinity and the second and third terms tend to zero. Thus there is an integer $K$ such that $b_{p(k)} < l_{k-1}$ for all $k \geq K$.

Let $0 < \delta < \alpha^{-l_K}$ and let $l(\delta)$ be the unique integer such that $\alpha^{-l(\delta) - 1} \leq \delta < \alpha^{-l(\delta)}$. Furthermore, let $k (\delta)$ be the unique integer such that $l_{k(\delta) - 1} < l (\delta) \leq l_{k (\delta)}$. Now, let $V$ be a component interval of $F_{a_{p(k(\delta))}}$ that is included in $V_{f^{-1} (k)}$. Then $V$ has diameter $\alpha^{-b(p(k(\delta))}$, and since $k (\delta) > K$ we have $b_{p(k(\delta))} < l(\delta) \leq l_{k(\delta)}$. Thus $V$ contains $2^{l(\delta) - b_{p(k(\delta))}}$ component intervals of $F_{a_{p(k(\delta))} + l (\delta) - b_{p(k(\delta))}}$, each of which has length $\alpha^{-l (\delta)}$ and contains points in $F$. Any interval of diameter $\delta$ can intersect at most two such intervals, so $N_{\delta} (F) \geq 2^{l(\delta) - b_{p(k(\delta))} - 1}$ and it follows that
\begin{equation*}
    \frac{\log N_{\delta} (F)}{-\log \delta} \geq \frac{(l(\delta) - b_{p(k(\delta))}-1) \log 2}{(l (\delta) + 1) \log \alpha} \geq \left( 1 - \frac{b_{p(k(\delta))}}{l_{k(\delta)-1}} - \frac{1}{l_{k(\delta)-1}} \right) \, \frac{\log2}{\log \alpha} \rightarrow \gamma
\end{equation*}
since $\frac{b_{p(k(\delta))}}{l_{k(\delta)-1}} \rightarrow 0$ by (\ref{limit}) and $\frac{\log 2}{\log \alpha} = \gamma$. Hence $\LBD F \geq \gamma$. \qed

\textit{Claim 3.} $\LMBD F = \gamma$.

\textit{Proof of Claim 3.} It suffices to show that $\LBD (F \cap U) = \gamma$ for all nonempty open subsets $U$ of $F$ (see \cite[Proposition~3.9]{falconer} for more details). Moreover, it follows by a version of Baire's Category Theorem that it is sufficient to show that $\LBD (F \cap V_w) = \gamma$ for all finite binary words $w$. To this end, let $w = w_1 \dots w_n \in \Sigma^{*}$ and define $g_w \colon \Sigma_2^{*} \rightarrow \Sigma_2^{*}$ by $g_w (v) = w v$. Now, for each word $u \in \Sigma_2^{*}$, define
\begin{equation*}
\begin{split}
    V^{'}_u &= V_{g_w (u)} \\
    a^{'}_u &= a_{g_w (u)} \\
    b^{'}_u &= b_{g_w (u)} \text{.}
\end{split}
\end{equation*}
Then applying the same arguments as in Claim 2, with $V_u$ replaced by $V^{'}_u$, $a_u$ replaced by $a^{'}_u$, and $b_u$ replaced by $b^{'}_u$ gives that $\LBD (F \cap V_w) = \gamma$.  Therefore we conclude that $\LMBD F = \gamma$. \qed

\textit{Claim 4.} $\HD F = 0$.

\textit{Proof of Claim 4.} We define a mass distribution on $F$ in an analogous manner to in the proof of Proposition \ref{C prop}. Namely, for each $n \in \mathbb{Z}^{+}$ and $w \in \Sigma_2^n$, let $\lambda_w$ denote the Lebesgue measure restricted to the interval $V_w$, normalised such that $\lambda_w (V_w) = 1$. Next, define the probability measure $\mu_n$ by $\mu_n = \frac{1}{2^n} \sum_{w \in \Sigma_2^n} \lambda_w$. Then there exists a probability measure $\mu$ such that $\mu_n$ converges weakly to $\mu$.

Let $x \in F$. Then there is an infinite binary word $w_1 w_2 \dots$ such that $x = \cap_n V_{w_1 \dots w_n}$. Let $(m_k)_k$ be the sequence of integers defined by $m_k= f^{-1} (w_1 \dots w_k)$ and for each $k \in \mathbb{Z}^{+}$ let $v_k (x)$ denote the interval $V_{w_1 \dots w_{a_{m_k}}}$. Each $v_k (x)$ has diameter $\alpha^{-b_{m_k}}$, therefore any ball with centre at $x$ and radius $r_k = \alpha^{b_{m_k}}$ can intersect at most three component intervals of $F_{a_{m_k}}$. Each has mass $2^{-a_{m_k}}$, so
\begin{equation*}
    \frac{\log \mu (B(x,r_k))}{\log r_k} \leq \frac{a_{m_k} \log 2 + \log 3}{b_{m_k} \log \alpha} \rightarrow 0 \text{.}
\end{equation*}
Hence it follows by the Mass Distribution Principle that $\HD F = 0$. \qed

\textit{Proof of Lemma \ref{F prop}.} Combining Claims 1, 3 and 4 completes the proof of Lemma \ref{F prop}. \qed

\subsection{Proof of Lemma \ref{LBD lemma}}

Finally, we present the proof of Lemma \ref{LBD lemma}. To ease notation, we write $C_n = C_n (r,t,\mathbf{n})$, $D_n = D (u,v,\mathbf{n},\mathbf{k})$, $E_n = E_n (w,\mathbf{j})$ and $F_n = F_n (s,\mathbf{a},\mathbf{b})$, for each integer $n$.

\textit{Proof of Lemma \ref{LBD lemma}.} Let $C$, $D$, $E$ and $F$ be as in the statement of Lemma \ref{LBD lemma}. In particular, recall that $\mathbf{n} = (n_k)_k$ is an increasing sequence of positive integers with $\frac{n_k}{n_{k+1}} \rightarrow 0$ as $k \rightarrow \infty$. It follows by monotonicity that $\LBD X \geq \LBD D = u$, so it remains to show the reverse inequality. To this end, for each integer $k$, let $\delta_k = 2^{-(n_{2k+1}-n_{2k})/u}$. Then, for all $\beta \leq u$ and $\gamma \leq v$, one interval of diameter $\delta_k$ will cover any component interval of $C_{n_{2k+1}} (\beta,\gamma,\mathbf{n})$. Hence $N_{\delta_k} (C) \leq 2^{n_{2k+1}}$, and since $D_{n_{2k+1}} (u,v,\mathbf{n},\mathbf{k}) \subseteq C_{n_{2k+1}} (u,v,\mathbf{n})$ and $F_{n_{2k+1}} (s,\mathbf{a},\mathbf{b}) \subseteq C (s,s,\mathbf{n})$ for all $k$, we have $N_{\delta_k} (D) \leq 2^{n_{2k+1}}$ and $N_{\delta_k} (F) \leq 2^{n_{2k+1}}$. For the bound on $N_{\delta_k} (E)$, observe that since $\UBD E = 0 < u$ there is an integer $K$ such that $N_{\delta_k} (E) \leq \delta_k^{-u} \leq 2^{n_{2k+1}}$ for all $k \geq K$. Hence, for $k \geq K$ we have $N_{\delta_k} (X) \leq N_{\delta_k} (C) + N_{\delta_k} (D) + N_{\delta_k} (E) + N_{\delta_k} (F) \leq 4 \cdot 2^{n_{2k+1}}$, so
\begin{equation*}
    \frac{\log N_{\delta_k} (X)}{-\log \delta_k} \leq u \, \left( \frac{n_{2k+1} + 2}{n_{2k+1} - n_{2k}} \right) \rightarrow u \text{,}
\end{equation*}
and we conclude that $\LBD X \leq u$. This completes the proof of Lemma \ref{LBD lemma}. \qed

\section{Outlook} \label{discussion}

\subsection{Higher dimensional analogues} \label{HD analogues}

One may ask whether an analogous result to Theorem \ref{main-thm} holds in higher dimensions. By modifying our construction, it is in fact possible to obtain the following result.
\begin{theorem} \label{HD-thm}
Let $r,s,t,u,v,w \in (0,d]$ with

\begin{center} 
{\begin{tikzpicture}
\node[left] at (0,0) {$r \leq s$};
\node at (0.1,0.4) [rotate = 40]{$\leq$};
\node at (0.1,-0.4)[rotate = -40]{$\leq$};
\node at (0.5,0.75) {$u$};
\node at (0.5,-0.75){$t$};
\node at (0.9,0.4) [rotate = -40]{$\leq$};
\node at (0.9,-0.4) [rotate = 40]{$ \leq$};
\node[right] at (1,0) {$v \leq w$.};
\end{tikzpicture}}
\end{center}
\noindent Then there exists a compact set $X \subseteq [0,1]^{d}$ with
\begin{equation*}
    \HDthm X = r \text{,}
\end{equation*}
\begin{equation*}
    \LMBDthm X = s \text{,}
\end{equation*}
\begin{equation*}
    \PDthm X = t \text{,}
\end{equation*}
\begin{equation*}
    \LBDthm X = u \text{,}
\end{equation*}
\begin{equation*}
    \UBDthm X = v \text{,}
\end{equation*}
\begin{equation*}
    \ADthm X = w \text{.}
\end{equation*}
Moreover, $X$ is an $r$-Hausdorff set and a $t$-packing set.
\end{theorem}
As in the proof of Theorem \ref{main-thm}, the set $X$ is constructed by taking the union of four sets $C$, $D$, $E$ and $F$, which could be thought of as higher dimensional analogues of the sets constructed in the proof of Theorem \ref{main-thm}. Instead of beginning the constructions of $C$ and $F$ with the unit interval we start with the unit $d$-cube, and for $D$ and $E$ we start with the corners of the cube, as opposed to the endpoints of the unit interval. We then construct the sets in an analogous manner, adapting the definitions of the contraction maps accordingly to take into account the extra dimensions. The proof follows similarly to that of Theorem \ref{main-thm}, but is somewhat more notationally awkward.

\subsection{Hewitt-Stromberg measures}

The so-called Hewitt-Stromberg measures provide a bridge between the Hausdorff and packing measures of a set. They were first introduced by Hewitt and Stromberg in their classical textbook \cite[(10.51)]{hewitt-stromberg} and have received attention in the fractal geometry community in recent years, for example in \cite{haase}, \cite{jurina}, \cite{zindulka}. For $\alpha > 0$ and $A \subseteq \mathbb{R}^d$ we denote the $\alpha$-dimensional lower and upper Hewitt-Stromberg measures of $A$ by $\mathcal{U}^{\alpha} (A)$ and $\mathcal{V}^{\alpha} (A)$ respectively. These measures satisfy $\mathcal{H}^{\alpha} (A) \leq \mathcal{U}^{\alpha} (A) \leq \mathcal{V}^{\alpha} (A) \leq \mathcal{P}^{\alpha} (A)$. The reader is referred to Edgar's book \cite{edgar} for a systematic introduction to the Hewitt-Stromberg measures.

The lower Hewitt-Stromberg measure can be used to give an alternative characteristaion of the lower modified box dimension, namely: $\LMBD A = \inf \{ \alpha \colon \mathcal{U}^{\alpha} (A) = 0 \} = \sup \{ \alpha \colon \mathcal{U}^{\alpha} (A) = \infty \}$. Therefore, it is natural to ask whether the set $X = C \cup D \cup E \cup F$ constructed in the proof of Theorem \ref{main-thm} satisfies $0 < \mathcal{U}^s (X) < \infty$. It turns out that this is not the case in general. It can be shown that $\mathcal{U}^s (X) < \infty$ through an identical argument to that used to show $\mathcal{P}^t (X) < \infty$; however, it can also be shown that $\mathcal{U}^s (F) = 0$. Since $\mathcal{U}^s (D) = \mathcal{U}^s (E) = 0$, it follows that $\mathcal{U}^s (X) > 0$ if and only if $\mathcal{U}^s (C) > 0$, and this is only true when $s = r$.

It is also known that the critical value of the upper Hewitt-Stromberg measure coincides with the packing dimension; however, the packing measure and the upper Hewitt-Stromberg measure do not coincide in general. Therefore one might ask whether $0 < \mathcal{V}^t (X) < \infty$. This is indeed the case. It follows immediately from our results that $\mathcal{V}^t (X) < \infty$, and it can be shown that $\mathcal{V}^t (C) > 0$. Hence $X$ is a $t$-upper Hewitt-Stromberg set, where, in an analogous manner to the Hausdorff and packing measures, a $t$-upper Hewitt-Stromberg set is one with positive and finite $t$-dimensional upper Hewitt-Stromberg measure.

\subsection{Homogeneity properties}

Finally, one could ask whether our set exhibits homogeneity properties similar to those satisfied by the set constructed by Nilsson and Wingren in \cite{nilsson-wingren}. In their paper, Nilsson and Wingren show that given any three numbers $r,s,t \in (0,d]$ with $r < s < t$, it is possible to construct a compact subset $K$ of $\mathbb{R}^d$ with $\HD (K \cap U) = r$, $\LBD (K \cap U) = s$ and $\UBD (K \cap U) = t$ for every open set $U$ with $K \cap U \neq \varnothing$. One can observe that such a result does not hold for our set. In fact, it is not even possible to obtain such a result, in general, for all the dimensions that we consider. It is known that if $\LBD (X \cap U) = \alpha$ for all open sets $U$ with $X \cap U \neq \varnothing$, then $\LMBD X = \LBD X = \alpha$ (there is an analogous result for the packing dimension and upper box dimension). Therefore, if $\LMBD X < \LBD X$ then it is impossible to have $\LMBD (X \cap U) = \LBD (X \cap U)$ for all open sets $U$ that intersect $X$ (and similarly for the packing dimension and upper box dimension).

\section*{Acknowledgements}

A. Mitchell would like to thank Jonathan Fraser and Kenneth Falconer for valuable discussions while writing this paper.

\end{spacing}

\end{document}